\documentclass[11pt, reqno]{amsart}
\usepackage{graphicx} 
\usepackage{amsmath,amssymb,amscd,mathrsfs,amscd}
\usepackage{graphics,verbatim, tikz}

\usepackage{hyperref}
\usepackage{pifont}

\usepackage{enumerate}

\usepackage{graphicx}
\usepackage{epstopdf}
\DeclareMathAlphabet{\mathpzc}{OT1}{pzc}{m}{it}
\usepackage{tikz,tikz-cd, color}
\usepackage[makeroom]{cancel}
\usetikzlibrary{arrows}

\usepackage{mathdots}

\usepackage{float}
\usepackage{caption, subcaption}
\usepackage{rotating}
\usepackage{pdfpages}
\allowdisplaybreaks

\usepackage{adjustbox}
\usetikzlibrary{matrix}

\tikzset{
  symbol/.style={
    draw=none,
    every to/.append style={
      edge node={node [sloped, allow upside down, auto=false]{$#1$}}}
  }
}

\usetikzlibrary{decorations.shapes}

\tikzset{decorate sep/.style 2 args=
{decorate,decoration={shape backgrounds,shape=circle,shape size=#1,shape sep=#2}}}

\newtheorem{theorem}{Theorem}[section]
\newtheorem{lemma}[theorem]{Lemma}
\newtheorem{proposition}[theorem]{Proposition}
\newtheorem{corollary}[theorem]{Corollary}

\newtheorem{definition}[theorem]{Definition}

\newtheorem{remark}[theorem]{Remark}

\theoremstyle{definition}

\newcommand{\bC}{\mathbb{C}}

\newcommand{\bR}{\mathbb{R}}

\newcommand{\cC}{\mathcal{C}}
\newcommand{\cD}{\mathcal{D}}

\newcommand{\cM}{\mathcal{M}}

\newcommand{\cP}{\mathcal{P}}

\newcommand{\cU}{\mathcal{U}}

\newcommand{\End}{\operatorname{End}}

\newcommand{\sgn}{\mathtt{sgn}}

\begin{document}
\title[On character table of Clifford groups]{On character table of Clifford groups}
\author[Lee, Yu, Peng, Lai]{Chin-Yen Lee, Wei-Hsuan Yu, Yung-Ning Peng, Ching-Jui Lai}
\address{Department of Mathematics, National Central University, Chung-Li, 32054, Taiwan.}
\email{109281001@cc.ncu.edu.tw (Lee), \quad whyu@math.ncu.edu.tw (Yu)}
\address{Department of Mathematics, National Cheng Kung University, Tainan, 70101, Taiwan. }
\email{ynp@gs.ncku.edu.tw (Peng), \quad  cjlai72@mail.ncku.edu.tw (Lai)}

\begin{abstract}
Based on a presentation of $\cC_n$ and the help of \cite{GAP}, we construct the character table of the Clifford group $\cC_n$ for $n=1,2,3$. As an application, we can efficiently decompose the (higher power of) tensor product of the matrix representation in those cases. Our results recover some known results in \cite{HWW, WF} and reveal some new phenomena. We prove that when $n\geq 3$, (1) the trivial character is the only linear character for $\cC_n$ and hence $\cC_n$ equals to its commutator subgroup, (2) the $n$-qubit Pauli group $\cP_n$ is the only proper non-trivial normal subgroup of $\cC_n$, (3) the matrix representation $\cM_{2^n}$ is a faithful representation for $\cC_n$. As a byproduct, we give a presentation of the finite symplectic group $Sp(2n,2)$ in terms of generators and relations.
\end{abstract}


\maketitle

\section{Introduction}

The Clifford group $\cC_n$ is one the most important algebraic structures in recent development of quantum information.
Not only the group itself, also its representation theory is playing an increasingly prominent role in quantum information.
Among the various representations of $\cC_n$, the matrix representation, or equivalently the tensor product of the standard representation and its dual, is known to have important real world applications in, for instances, quantum computing \cite{CTV, F+7}, quantum error correction \cite{Go1, Go2},
quantum image process \cite{McBS+,LC}, 
randomized benchmarking \cite{HWFW,HWW,WF}, 
and particularly unitary $t$-design \cite{BOZ,NWD,We,Zh,ZKGG}.

As a consequence, to have a better understanding about the group structure of $\cC_n$ and its representation theory are natural and urgent tasks for us. Indeed, the structure and properties of $\cC_n$ have been studied in \cite{Ap,GS,Gr,KMP,Oz,Sel} from different points of view with various emphases and applications. 
On the other hand, the matrix representation $\cM_{2^n}$ of $\cC_n$ are intensively studied for they are extremely useful in many applications mentioned above. In particular, the decomposition of $\cM_{2^n}$ into irreducible $\cC_n$-invariant subspaces is studied in \cite{WF}, the decomposition of the two-copy tensor product $(\cM_{2^n})^{\otimes 2}$ is studied in \cite{HWW}, see also \cite{SH} for some partial results and conjectures about its three-copy tensor product $(\cM_{2^n})^{\otimes 3}$.

Since $\cC_n$ is a finite group over the field of complex numbers, its character table is the most essential information in the study of its finite dimensional representations. It is a surprising fact that their character tables are still unavailable in literature, except for $n=1,2$ that were constructed or observed very recently \cite{KMP,Ma}.

In this article, we construct the character table of $\cC_n$ for $n\leq 3$ and 
as a consequence, we are able to decompose any finite dimensional representation into a direct sum of its irreducible constituents simply by solving linear equations.
As illustrating examples, we decompose $(\cM_{2^n})^{\otimes m}$ for several $m$ for each $1\leq n\leq 3$. 
For $m=1,2$, our result recovers the results in \cite{HWW,WF}, and our approach works uniformly well for any reasonably large $m\geq 3$.

Our main tool is a concrete presentation of $\cC_n$, see Theorem \ref{presc}, together with the help of \cite{GAP}.
The case for $n=1$ is actually well known, for $\cC_1$ is isomorphic to the symmetric group $S_4$ \cite{KMP}. 
For $n=2$, the character table was independently obtained in a very recent article \cite{Ma} from a very different point of view; see Remark \ref{cntab}. 
We should emphasize here that the properties of the symplectic group $Sp(2n,2)$ is fundamental in almost all previous studies for $\cC_n$ and its representations. 
In particular, based on the structure of finite symplectic groups, there have been a series of interesting works studying the commutant algebra of higher tensor product of the matrix representation (more precisely, the Schur-Weyl duality; see \cite{FH} for classical cases) of $\cC_n$, together with the mixed version of Schur-Weyl duality and their applications; see \cite{GNW,MMG1,MMG2}.
In this article we provide an alternative approach to study the representation theory of $\cC_n$, which does not rely on any deep results about the symplectic groups. As a byproduct, based on the fact that the $Sp(2n,2)$ is a quotient of $\cC_n$, we can modify the presentation given in Theorem~\ref{presc} to obtain a presentation of $Sp(2n,2)$ in terms of abstract generators and relations, see Theorem~\ref{presp}, which is of independent interest and might be useful elsewhere.

We should emphasize that the case $n=3$ reveals some features that do not appear in the previous small rank cases. 
In fact, some features in the cases $n=1,2$ are no longer true in the rest members of the family, 
where some general features start to appear when $n=3$, see Theorem \ref{c=cc}. 
Based on observations of our results, we can further prove some general facts about $\cC_n$, see Theorem \ref{normalc} and Corollary \ref{faith}. 
We expect that a better understanding of the Clifford groups and their representations will potentially improve the performances, 
such as better efficiency or accuracy, in their various practical applications.

The article is organized as follows.
We setup the notation and recall the definition of the Clifford groups in section 2. 
The very first case $\cC_1$ is explained in section 3, together with a decomposition of the tensor products of the matrix representation.
In section 4 we provide a presentation for general $\cC_n$. 
In section 5 we use the presentation to obtain the character table for $\cC_2$ and hence a decomposition of the tensor products of the matrix representation as well. The case for $\cC_3$ is similarly obtained, and then some facts about $\cC_n$ are proved for all $n\geq 3$. 
A presentation of $Sp(2n,2)$ is provided in section 6.

\section*{Acknowledgments}
We thank Prof. Ming-Hsiu Hsieh, the director of the Hon Hai Quantum Computing Research Center in Taiwan, for bringing our attention to the theory of unitary $t$-design. 
We thank Jeng-Daw Yu for technical support.
Yu is partially supported by the NSTC (Grant No. 109-2628-M-008-002-MY4); 
Peng is partially supported by the NSTC (Grant No. 111-2628-M-008-001-MY3) and the National Center of Theoretical Sciences, Taipei, Taiwan;
Lai is partially supported by the NSTC (Grant No. 111-2115-M-006-002-MY2) and QFort at National Cheng Kung University.
We are grateful to Jason Saied and Markus Heinrich for their comments after we posted the first draft of the article on the arXiv, that turn our previous conjectures into solid affirmative theorems; see Remark \ref{conj}. We also thank Shih-Chang Huang and Shu-Yen Pan for helping us find a counterexample of one previous conjecture; see Remark \ref{intab}.

\section{Preliminaries}
In this section we set our notation and recall some basic facts of Clifford groups. We start with the definition of Pauli matrices.

\begin{definition}
The single-qubit Pauli matrices consist of 2 by 2 complex matrices 
\[
\hat\cP_1=\{ I_2= \left(\begin{array}{cc}
    1 & 0 \\ 
    0 & 1 \\ 
  \end{array}\right), \,\,
X=\left(\begin{array}{cc}
    0 & 1 \\ 
    1 & 0 \\ 
  \end{array}\right),\,\,
Y=\left(\begin{array}{cc}
    0 & -i \\ 
    i & 0 \\ 
  \end{array}\right),\,\,
Z=\left(\begin{array}{cc}
    1 & 0 \\ 
    0 & -1 \\ 
  \end{array}\right)\}
\]
The $n$-qubit Pauli matrices is the set $\hat\cP_n:=\{\sigma_1\otimes \cdots \otimes \sigma_n \, | \, \sigma_i \in \cP_1\}\subseteq (\cU_2(\bC))^{\otimes n}$. 
\end{definition}

By matrix multiplication, it is easy to see that $\hat\cP_1$ generates a subgroup of order 16 in $\cU_2(\bC)$, the group of 2 by 2 unitary matrices. In particular, the center of this subgroup is precisely $U_1=\{c \cdot I_2 \, | \, c\in\{\pm1, \pm i\} \}$, called the global phase. 
In real world applications, the global phase is usually ignored. The quotient group $\cP_1:=\langle \hat\cP_1 \rangle/ U_1$ is called the {\em single-qubit Pauli group}, which is isomorphic to the Klein 4 group. Similarly, the {\em n-qubit Pauli group} is defined by $\cP_n:=\langle \hat\cP_n\rangle/U_1$, where we implicitly identify $U_1$ with the center of $\langle \hat\cP_n\rangle$, which consists of the scalar matrices $\{c \cdot (I_2)^{\otimes n} \, | \, c\in\{\pm1, \pm i\} \}$. We have $\cP_n\cong (\cP_1)^{\oplus n}$, an abelian group.

\begin{definition}
The $n$-qubit Clifford group $\cC_n$ is the normalizer of $\cP_n$ (ignoring the phase). That is,
\[
\cC_n=\{ A\in (\cU_2(\bC))^{\otimes n}\, | \, A\sigma A^{-1} \in \langle \hat\cP_n \rangle, \text{ for all } \sigma \in \langle \hat\cP_n \rangle \}/ 
\langle \omega \cdot (I_2)^{\otimes n}\rangle, 
\]
where $\omega=e^{\frac{\pi i}{4}}$ is a primitive $8^{th}$ root of 1.
\end{definition}

The order of $\cC_n$ is well known, see for example \cite{Oz} for a detailed proof.
\begin{equation}\label{ocn}
|\cC_n|=2^{n^2+2n}\prod_{j=1}^n (2^{2j}-1)
\end{equation}
We list the order of the first few members of this family in the following table.
\begin{center}
\begin{tabular}{|c|c|c|c|c|c|}
\hline
$n$ & 1 & 2 & 3 & 4 & 5  \\
\hline
 $|\cC_n|$ & 24 & 11520 & 92897280 & 12128668876800 & 25410822678459187200 \\
\hline
\end{tabular}\\[.5cm]
\end{center}

\section{The single-qubit case}\label{C1}

It is well known \cite{GS, Oz} that $\cC_1$ can be described as the group of rotations of the Bloch sphere in $\bR^3$ that permute $\pm x$, $\pm y$, and $\pm z$ directions, but it is less noticed that $\cC_1$ is isomorphic to $S_4$ \cite{KMP}, the symmetric group of 4 objects. Later we will give a presentation for general $\cC_n$ and deduce some further results. Before that we first provide some detail in this case.

Consider the following two complex matrices, called the Hadamard gate and the phase gate, respectively:
\begin{equation}\label{HP}
H=
\frac{1}{\sqrt{2}}\left(
\begin{array}{cc}
  1 &  1     \\
  1  & -1    
\end{array}
\right),
\qquad
P=
\left(
\begin{array}{cc}
  1 &  0     \\
  0  & i    
\end{array}
\right).
\end{equation}
It can be checked by definition that both of $H$ and $P$ are elements in the single-qubit Clifford group.\footnote{The matrix equalities in this article are mostly in $\cC_n$, a quotient group, which means that two matrices are considered equal if they differ by a scalar matrix. For example, $XY=iZ$ as complex matrices but we have $XY=Z$ as elements in $\cC_1$. It shall be clear from the context that we are dealing with the whole unitary group or its quotient so we use the same notation for a unitary matrix and its image after ignoring the global phase.} In fact, these two elements are sufficient to generate the whole $\cC_1$. 
\begin{lemma}
A presentation of $S_4$ is given by
\begin{center}
$\langle t,s \, | \, t^2=s^4=1, \, (ts)^3=1, \, (tsts^{-1})^3=1, \, (ts^2ts^2)^2=1\rangle$ 
\end{center}
\end{lemma}
\begin{corollary}
$\cC_1\cong S_4$
\end{corollary}
\begin{proof}
The assignment $H\mapsto t$ and $P\mapsto s$ gives a group isomorphism.
\end{proof}

Since the character table of $S_4$ is well known, we can use it to decompose the tensor product of any finite dimensional representations of $\cC_1$.
In particular, we consider the matrix representation of $\cC_1$, denoted by $\cM_2$, and its tensor products $(\cM_2)^{\otimes m}$, for any $m\in\mathbb{N}$. 
Here $\cM_2$ is the set of all 2 by 2 complex matrices and $\cC_1$ acts on $\cM_2$ by conjugation:
\begin{equation}\label{conm2}
A \cdot M := A M A^{-1}, \qquad \forall A\in \cC_1, M\in \cM_2,
\end{equation}
and then $\cC_1$ acts on $\cM_2^{\otimes m}$ by iteration of the comultiplication. One should note that the global phase acts trivially in (\ref{conm2}) and hence the action indeed factors though the quotient.

It is well known that the number of irreducible characters is exactly the number of conjugacy classes of the group. For $S_4$, there are precisely 5 classes and let $\{\chi_i, 1\leq i\leq 5\}$ denote the irreducible characters of $S_4$. We use a non-standard labeling for the characters so we list their character values on each conjugacy class explicitly.
\begin{center}
\begin{tabular}{|c|c|c|c|c|c|}
\hline
Character $\backslash$ Representative & $e$ & $H$  &  $P^2$ & $HP$ & $P$ \\
\hline
 $\chi_1$ & 1 & 1 & 1 & 1 & 1 \\
\hline
 $\chi_2$ & 3 & 1 & -1  & 0 & -1\\
\hline
  $\chi_3$ & 2 & 0 & 2  & -1 & 0\\
\hline
 $\chi_4$ & 3 & -1 & -1  & 0 & 1\\
\hline
 $\chi_5$ & 1 & -1 & 1  & 1 & -1\\
\hline
\end{tabular}\\[.5cm]
\end{center} 
One can compute the character value of $\text{ch}(\cM_2)$ by definition of the action given in (\ref{conm2}):
\begin{center}
\begin{tabular}{|c|c|c|c|c|c|}
\hline
   &  $e$ & $H$  &  $P^2$ & $HP$ & $P$  \\
\hline
 $\text{ch}(\cM_2)$ & 4 & 0 & 0 & 1 & 2 \\
\hline
\end{tabular}\\[.5cm]
\end{center} 

We use the following decomposition vector $v_m:=(a_m,b_m,c_m,d_m,e_m)\in(\mathbb{Z}_{\geq 0})^5$ to denote the decomposition of the character into its irreducible constituents:
\[ \text{ch}\big((\cM_2)^{\otimes m}\big)=a_m\chi_1+b_m\chi_2+c_m\chi_3+d_m\chi_4+e_m\chi_5.\] 
Using the rows of our character table as an orthonormal basis for the vector space of class functions on $S_4\cong \cC_1$, one can easily obtain the decomposition of $(\cM_2)^{\otimes m}$ for any given $m$ simply by solving linear equations rather than finding $\cC_1$-invariant subspaces in $(\cM_2)^{\otimes m}$. We list the first several decomposition vectors as follows:
\begin{align*}
v_1&=(1,0,0,1,0)\\
v_2&=(2,1,1,3,0)\\
v_3&=(5,6,5,10,1)\\
v_4&=(15,28,21,36,7)\\
v_5&=(51,120,85,136,35)\\
v_6&=(187,496,341,528,155)
\end{align*}
In particular, our $v_1$ and $v_2$ recover the results in \cite{WF} and \cite[Theorem 1]{HWW} for single-qubit Clifford group, respectively.
\begin{remark}
By computation, we found a recursive formula of $v_m$ as follows:
\[
v_{m+1}=(a_m+d_m,2b_m+c_m+d_m+e_m,b_m+c_m+d_m,a_m+b_m+c_m+2d_m,b_m+e_m)
\]
The formula looks symmetric using our labeling of irreducible characters, but due to the non-symmetric initial vector $v_1$, the resulted vectors are not quite symmetric.
\end{remark}

\section{A presentation for Clifford groups}
The following result provides a foundation of this article.
\begin{theorem}\cite[Proposition 7.1]{Sel}\label{presc}
The group $\mathcal{C}_n$ is generated by the symbols 
$$\{H_i, P_i, Z_{j} \,\, | \,\, 1\leq i\leq n, \, 1\leq j\leq n-1\}$$ subjected to the following relations:
\begin{itemize}
\item(R1) $Z_j^2=H_i^2=P_i^4=1$ for all $1\leq i\leq n$ and $1\leq j\leq n-1$. \\[-3mm]
\item(R2) $(H_iP_i)^3=1$ for all $i$.  \\[-3mm]
\item(R3) $(H_iP_iH_iP_i^3)^3=1$ for all $i$. \\[-3mm]
\item(R4) $(H_iP_i^2H_iP_i^2)^2=1$ for all $i$. \\[-3mm]
\item(R5) $H_iH_j=H_jH_i$, \,\, $P_iP_j=P_jP_i$, \,\, and \,\, $Z_iZ_j=Z_jZ_i$ for all $i\neq j$. \\[-3mm]
\item(R6) $Z_jP_i=P_iZ_j$ for all $i, j$. \\[-3mm]
\item(R7) $Z_jH_{i}=H_{i}Z_j$ for all $i\notin\{j,j+1\}$.\\[-3mm]
\item(R8) $Z_jH_jP_j^2H_j=H_jP_j^2P_{j+1}^2H_jZ_j$ for all $1\leq j\leq n-1$.\\[-3mm]
\item(R9) $Z_jH_{j+1}P_{j+1}^2H_{j+1}=H_{j+1}P_j^2P_{j+1}^2H_{j+1}Z_j$ for all $1\leq j\leq n-1$.\\[-3mm]
\item(R10) $Z_jH_jZ_j =P_jH_jP_jP_{j+1}Z_jH_jP_j$  for all $1\leq j\leq n-1$.\\[-3mm]
\item(R11) $Z_jH_{j+1}Z_j =P_{j+1}H_{j+1}P_jP_{j+1}Z_jH_{j+1}P_{j+1}$ for all $1\leq j\leq n-1$.\\[3mm]
When $n\geq 3$, we also have
\item(B1) $Z_jH_jH_{j+1}Z_jH_{j+1}H_{j+2}Z_{j+1}H_{j+1}H_{j+2}Z_jH_jH_{j+1}Z_j=\\[1mm]
		  Z_{j+1}H_{j+1}H_{j+2}Z_{j+1}H_jH_{j+1}Z_jH_jH_{j+1}Z_{j+1}H_{j+1}H_{j+2}Z_{j+1}$ for all $1\leq j\leq n-2$.\\[-3mm]
\item(B2) $(Z_{j+1}H_jH_{j+1}Z_jH_jH_{j+1}Z_j)^3=1$ for all $1\leq j\leq n-2$.\\[-3mm]
\item(B3) $(Z_jH_{j+1}H_{j+2}Z_{j+1}H_{j+1}H_{j+2}Z_{j+1})^3=1$ for all $1\leq j\leq n-2$.
\end{itemize}
\end{theorem}

\begin{remark}\label{s4com}
For each $1\leq k\leq n$, we have $T_k:=\langle H_k,P_k\rangle\cong \cC_1$. Moreover, if $k\neq j$ then the two subgroups $T_k$ and $T_j$ commute each other in $\cC_n$. 
\end{remark}
\begin{remark}\label{8times}
The Clifford groups considered in \cite{Sel} include the global phases as group elements and hence the order is 8 times of the group in our consideration. Our presentation here is modified from that given in \cite[Proposition 7.1]{Sel} by quotient out those scalars. 
Similar or equivalent presentations (with or without global phases) are used in a recent article \cite{KMP}, which focuses on the study of all subgroups of $\cC_2$, except that the $CNOT$ gate is used there instead of the controlled-$Z$ gate here. 
\end{remark}

\begin{proof}
Denote the abstract group defined by the presentation in Theorem \ref{presc} by $\cD_n$. One checks that the assignment $\cD_n \rightarrow \cC_n$ given by
\begin{align*}
H_i\mapsto& I_2^{\otimes (i-1)}\otimes  H \otimes  I_2^{\otimes (n-i)},\\
P_i\mapsto& I_2^{\otimes (i-1)}\otimes  P \otimes  I_2^{\otimes (n-i)},\\ 
Z_j\mapsto&  I_2^{\otimes (j-1)}\otimes CZ \otimes I_2^{\otimes (n-j-1)},
\end{align*}
where \[CZ= \left(\begin{array}{cc}
    1 & 0 \\ 
    0 & 0 \\ 
  \end{array}\right) \otimes I_2 +
\left(\begin{array}{cc}
    0 & 0 \\ 
    0 & 1 \\ 
  \end{array}\right) \otimes \left(\begin{array}{cc}
    1 & 0 \\ 
    0 & -1 \\ 
  \end{array}\right),\]
preserves the defining relations and hence we have a group epimorphism $\Phi:\cD_n\rightarrow \cC_n$. Here the matrices $H$ and $P$ are defined in \textsection \ref{C1}, and $CZ$ is known as the controlled-$Z$ gate. 
The injectivity of $\Phi$ follows from \cite[Proposition~7.1]{Sel} and the fact that the order of $\cD_n$ is $\frac{1}{8}$ of the number given in \cite[Corollary 5.6]{Sel}, which is exactly the order of $\cC_n$ given in (\ref{ocn}); see Remark \ref{8times}.
\end{proof}

\begin{remark}\label{swap}
Note that $(Z_jH_{j+1}H_j)^3\in \End \big( (\bC^2)^{\otimes n}\big)\cong (\End(\bC^2))^{\otimes n}$ is the operator swapping the $j^{th}$ and the $(j+1)^{th}$ tensor factor. As a consequence, all $H_i$'s are conjugate in $\cC_n$; similarly, all $P_i$'s are conjugate in $\cC_n$ and all $Z_j$'s are as well. Moreover, the wreath product $S_4 \wr S_n$ is naturally a subgroup of $\cC_n$.   
\end{remark}

\section{The cases for Two-qubit, Three-qubit, and beyond}
Note that the relations $(B1)$--$(B3)$ do not appear for $n=1,2$, which makes them different from other members in this family.
We have seen the complete character table for $\cC_1$ and now we analyze the case for $\cC_2$.

\subsection{Two-qubit case}
$\cC_2$ is generated by $\{H_i,P_i,Z=Z_1 \,|\, i=1,2\}$.
With the presentation given in Theorem \ref{presc} and the help of GAP, we find that there are only two non-trivial proper normal subgroups in $\cC_2$ with order 16 and 5760, denoted by $N_1$ and $N_2$, respectively.

The order of $N_1$ is 16, and by definition the Pauli group $\cP_2$ is a normal subgroup of $\cC_2$ with the same order, therefore they coincide.
It then follows that $N_1$ is generated by $\{X_i,Y_i \, | \, i=1,2\}$, where $X_1=X\otimes I_2$, $X_2=I_2\otimes X$ and similarly for $Y_i$.
We do not further investigate the structure of $\cC_2/N_1$ here, except pointing out that it is isomorphic to the symplectic group $Sp(4,\mathbb{Z}_2)$ \cite{Ap,BRW,Gr,GNW}, which is non-abelian, simple and of order 720; see also a recent article \cite{KMP} for a further discussion about all subgroups of $\cC_2$.

The order of $N_2$ is 5760, which is exactly half of the order of $\cC_2$. It in turns implies that $N_2$ is precisely the commutator subgroup $[\cC_2,\cC_2]$ since it is the only proper normal subgroup with abelian quotient. We collect these facts and their consequences as the following theorem. 

\begin{theorem}\label{c2thm}
The following facts hold in $\cC_2$.
\begin{enumerate}
\item There are only two proper nontrivial normal subgroups in $\cC_2$. One of them is the 2-qubit Pauli group, which has order 16. The other one is its commutator subgroup $[\cC_2,\cC_2]$, which has order 5760. 
\item There are 21 conjugacy classes in $\cC_2$ and hence 21 irreducible characters.
\item There are only two linear (1 dimensional) representations of $\cC_2$, the trivial representation and the sign representation $\sgn:\cC_2\rightarrow {\{\pm1\}}$ defined by $\sgn(H_i)=\sgn(P_i)=\sgn(Z)=-1$ for $i=1,2$.
\end{enumerate}
\end{theorem}
\begin{proof}
The first two statements follow from the calculation by GAP, where the last one follows from the fact that $|\cC_2/[\cC_2,\cC_2]|=2$ and the defining relations are preserved by the assignment of $\sgn$.
\end{proof}

We also obtain the character table of $\cC_2$ by GAP; see Table \ref{TabCC2}. To further explain it, we provide a table to label each conjugacy class and pick a representative element in it; see Table \ref{TabC2}.
\begin{table}[ht]
\caption{Conjugacy classes of $\cC_2$}
\begin{center}
\begin{tabular}{|c|c|c|c|c|c|}
\hline
 Label & Size & Representative  &  Label & Size & Representative \\
\hline
 No.1 & 1 & $e$ &  No.12 & 120 &   $H_1$    \\
\hline
 No.2 &  640 & $P_2P_1H_2H_1$  & No.13 & 160 &   $P_1H_1$ \\
\hline
  No.3 & 60 & $Z$ &   No.14 & 480 &  $P_2^2P_1H_1$    \\
\hline
  No.4 & 1920  & $ZH_2H_1$ &   No.15 & 960 &  $P_2P_1H_1$ \\
\hline
  No.5 & 15 & $P_1^2$  &  No.16 & 2304 &  $ZP_1H_2H_1$     \\
\hline
  No.6 & 180 & $P_2P_1$  &   No.17 & 180 &  $H_2H_1$ \\
\hline
  No.7 & 720 &  $ZH_1$  &   No.18 & 720 & $H_1H_2P_2P_1^{-1}ZH_1H_2Z$      \\
\hline
 No.8 & 30 & $P_1$  &   No.19 & 960 &   $P_1H_2H_1$   \\
\hline
 No.9 & 360 &   $P_2H_1$ &  No.20 & 720 &  $H_1ZH_2P_2Z$   \\
\hline
  No.10 & 180  & $H_2P_2^2H_2Z$  & No.21 & 720 &    $ZP_2H_1$   \\
\hline
  No.11 & 90  &   $P_2P_1^2$ &  & &  \\
\hline
\end{tabular}\\[1cm]
\end{center}
\label{TabC2}
\end{table}

\begin{table}[ht]
\caption{Character Table of $\cC_2$}
\begin{center}
\setlength\tabcolsep{3pt}
\begin{tabular}{|c|c|c|c|c|c|c|c|c|c|c|c|c|c|c|c|c|c|c|c|c|c|}
\hline
  Character $\backslash$ Class  & 1 & 2&3&4&5&6&7&8&9&10&11&12&13&14&15&16&17&18&19&20&21 \\
\hline\hline
 $\chi_1$ & 1 & 1&1&1&1&1&1&1&1&1&1&1&1&1&1&1&1&1&1&1&1\\
\hline
$\chi_2$ & 1 & 1&-1&-1&1&1&1&-1&1&-1&-1&-1&1&1&-1&1&1&1&-1&-1&-1\\
\hline
$\chi_3$ & 5 & -1&-1&-1&5&1&-1&3&1&-1&3&3&2&2&0&0&1&-1&0&1&1\\
\hline
$\chi_4$ & 5& -1& 1& 1& 5& 1& -1& -3& 1& 1& -3& -3& 2& 2& 0& 0& 1& -1& 0& -1& -1\\
\hline
$\chi_5$ & 5 & 2 & -3 & 0 & 5 & 1 & -1 & 1 & 1 & -3 & 1 & 1 & -1 & -1 & 1 & 0 & 1 & -1 & 1 & -1 & -1\\
\hline
$\chi_6$ & 5 & 2 & 3 & 0 & 5 & 1 & -1 & -1 & 1 & 3 & -1 & -1 & -1 & -1 & -1 & 0 & 1 & -1 & -1 & 1 & 1\\
\hline
 $\chi_7$ & 9 & 0 & -3 & 0 & 9 & 1 & 1 & -3 & 1 & -3 & -3 & -3 & 0 & 0 & 0 & -1 & 1 & 1 & 0 & 1 & 1\\
\hline
$\chi_8$ & 9 & 0 & 3 & 0 & 9 & 1 & 1 & 3 & 1 & 3 & 3 & 3 & 0 & 0 & 0 & -1 & 1 & 1 & 0 & -1 & -1 \\
\hline
$\chi_9$ &10 & 1 & -2 & 1 & 10 & -2 & 0 & 2 & -2 & -2 & 2 & 2 & 1 & 1 & -1 & 0 & -2 & 0 & -1 & 0 & 0 \\
\hline
$\chi_{10}$ & 10 & 1 & 2 & -1 & 10 & -2 & 0 & -2 & -2 & 2 & -2 & -2 & 1 & 1 & 1 & 0 & -2 & 0 & 1 & 0 & 0\\
\hline
$\chi_{11}$ & 15 & 0 & 3 & 0 & -1 & -1 & -1 & -5 & -1 & -1 & 3 & -1 & 3 & -1 & 1 & 0 & 3 & 1 & -1 & 1 & -1 \\
\hline
$\chi_{12}$ & 15 & 0 & -3 & 0 & -1 & 3 & 1 & -7 & -1 & 1 & 1 & 1 & 3 & -1 & -1 & 0 & -1 & -1 & 1 & 1 & -1 \\
\hline
$\chi_{13}$ & 15 & 0 & -3 & 0 & -1 & -1 & -1 & 5 & -1 & 1 & -3 & 1 & 3 & -1 & -1 & 0 & 3 & 1 & 1 & -1 & 1 \\
\hline
$\chi_{14}$ & 15 & 0 & 3 & 0 & -1 & 3 & 1 & 7 & -1 & -1 & -1 & -1 & 3 & -1 & 1 & 0 & -1 & -1 & -1 & -1 & 1 \\
\hline
$\chi_{15}$ & 16 & -2 & 0 & 0 & 16 & 0 & 0 & 0 & 0 & 0 & 0 & 0 & -2 & -2 & 0 & 1 & 0 & 0 & 0 & 0 & 0 \\
\hline
$\chi_{16}$ & 30 & 0 & -6 & 0 & -2 & 2 & 0 & -2 & -2 & 2 & -2 & 2 & -3 & 1 & 1 & 0 & 2 & 0 & -1 & 0 & 0 \\
\hline
$\chi_{17}$ & 30 & 0 & 6 & 0 & -2 & 2 & 0 & 2 & -2 & -2 & 2 & -2 & -3 & 1 & -1 & 0 & 2 & 0 & 1 & 0 & 0 \\
\hline
$\chi_{18}$ & 45 & 0 & -3 & 0 & -3 & 1 & -1 & 9 & 1 & 1 & 1 & -3 & 0 & 0 & 0 & 0 & -3 & 1 & 0 & 1 & -1  \\
\hline
$\chi_{19}$ & 45 & 0 & 3 & 0 & -3 & 1 & -1 & -9 & 1 & -1 & -1 & 3 & 0 & 0 & 0 & 0 & -3 & 1 & 0 & -1 & 1\\
\hline
$\chi_{20}$ & 45 & 0 & -3 & 0 & -3 & -3 & 1 & -3 & 1 & 1 & 5 & -3 & 0 & 0 & 0 & 0 & 1 & -1 & 0 & -1 & 1  \\
\hline
$\chi_{21}$ & 45 & 0 & 3 & 0 & -3 & -3 & 1 & 3 & 1 & -1 & -5 & 3 & 0 & 0 & 0 & 0 & 1 & -1 & 0 & 1 & -1\\
\hline
\end{tabular}\\[.5cm]
\end{center} 
\label{TabCC2}
\end{table}

With the character table, we may decompose any finite dimensional $\cC_2$-module $V$ into a direct sum of irreducible submodules by solving linear equations, as long as the character of $V$ is known. Here again we consider the tensor product of the matrix representation of $\cC_2$, denoted by $\cM_{4}$. As a set, $\cM_4=\cM_2\otimes \cM_2$, and $\cC_2$ acts on $\cM_4$ by conjugation on each tensor factor:
\[
(A_1\otimes A_2) \cdot (M_1\otimes M_2) := (A_1 M_1 {A_1}^{-1}) \otimes (A_2 M_2 A_2^{-1}). 
\]
As before, $\cC_2$ acts on $(\cM_4)^{\otimes m}$ by iteration of the comultiplication. Let $\chi=\text{ch}(\cM_4)$ denote the character of $\cM_4$. 
We compute the character values of $\chi$ on each classes by definition, using the representatives given in Table \ref{TabC2}, and it turns out that $\chi=\chi_1+\chi_{14}$. We further compute the decomposition of $\chi^{\otimes m}$ for $m=2,3,4,5$, and the results are listed below. 
\begin{align*}
v_1=&(1,0,0,0,0,0,0,0,0,0,0,0,0,1,0,0,0,0,0,0,0)\\
v_2=&(2,0,1,0,0,0,0,1,0,0,0,0,1,4,0,0,1,2,0,0,1)\\
v_3=&(6,0,8,0,2,1,0,9,4,1,1,0,12,23,4,6,14,27,6,9,18)\\
v_4=&(29 ,  0,  70 ,  7 , 35,  21 , 21 , 91,  70,  35,  42,  28, 140, 190 , 84 ,140 ,196 ,350 ,168, 210, 280)\\
v_5=&( 219 ,  28 , 750 , 245 , 525 , 385 , 567, 1107 ,1050 , 735 , 980,  840, 1800, 2076, 1428, 2520 ,\\
&\,2920, 4860 ,3360 ,3780, 4320)
\end{align*}
Again, our $v_1$ and $v_2$ recover the results in \cite{WF} and \cite[Theorem 1]{HWW}\footnote{We found a typo in Appendix B of \cite{HWW}. The dimension of $V_{\{\text{adj}\}}^\perp$ should be half of the number provided there.}  for two-qubit Clifford group, respectively.

\subsection{Three-qubit case and beyond}
With the presentation given in Theorem \ref{presc}, we use GAP to analyze the structure of $\cC_3$.
As a result, there are 67 conjugacy classes in $\cC_3$ and there is only one proper normal subgroup of order 64, which must be the 3-qubit Pauli group $\cP_3$. 
We also obtain the character table for $\cC_3$ through GAP; see \textsection\ref{appc3}.
In particular, the second conjugacy class (which contains 63 elements) in the table, together with the identity class, is precisely the normal subgroup $\cP_3$.

With the character table, we are able to decompose any finite dimensional $\cC_3$-modules as long as we know their character values.
Again we consider the matrix representation $\cM_8=(\cM_2)^{\otimes 3}$, where $\cC_3$ acts on it by conjugation on each tensor factor just like the previous cases and similarly for $(\cM_8)^{\otimes m}$. 
The decomposition of $(\cM_8)^{\otimes m}$ for any reasonably large $m$ can be obtained by solving linear equations in 67 variables.
For $m=1,2,3$ the results are listed below as illustrating examples.
\begin{align*}
v_1=&\,\, \chi_1+\chi_{10},\\
v_2=&\,\, 2\chi_1+\chi_6+\chi_7+4\chi_{10}+\chi_{28}+\chi_{35}+\chi_{40}+\chi_{44}+\chi_{47},\\
v_3=&\,\, 6\chi_1+\chi_5+9\chi_6+9\chi_7+24\chi_{10}+\chi_{13}+2\chi_{15}+\chi_{16}+4\chi_{17}+4\chi_{18}+4\chi_{22}+4\chi_{25}\\
        & +17\chi_{28}+\chi_{30}+2\chi_{31}+18\chi_{35}+18\chi_{40}\chi_{43}+18\chi_{44}+6\chi_{45}+4\chi_{46}+36\chi_{47}+9\chi_{50}\\ 
        & +4\chi_{51}+12\chi_{53}+6\chi_{57}+9\chi_{58}+12\chi_{60}+6\chi_{63}+18\chi_{65}+9\chi_{66}.
\end{align*}

We point out here that, different from previous $\cC_1$ and $\cC_2$, the trivial character is the only linear character for $\cC_3$. In particular, the relations $(B2)$-$(B3)$ are not preserved by the sign representation. From this observation, one can in fact prove the following general result.
\begin{theorem}\label{c=cc}
For $n\geq 3$, the trivial character is the only linear character of $\cC_n$ and hence $\cC_n=[\cC_n,\cC_n]$.
\end{theorem}
\begin{proof}
A linear character $\chi:\cC_n\rightarrow \bC$ is uniquely determined by its values on the generating set $\{H_k,P_k,Z_j  \,|\, 1\leq k \leq n, 1\leq j\leq n-1\}$.
Note that $\chi(H_k)$ and $\chi(Z_j)$ can only be $\pm1$, while $\chi(P_k)\in\{\pm1,\pm i\}$. 

Suppose that $\chi(P_k)\neq 1$ for some $1\leq k\leq n$. 
Consider the subgroup $T_k$ of $\cC_n$ generated by $H_k, P_k$. This subgroup is isomorphic to $\cC_1$, and the restriction $\chi|_{T_k}$ is a non-trivial linear character of $\cC_1$, which must be the sign character according to our discussion in \textsection\ref{C1}, forcing $\chi(P_k)=-1=\chi(H_k)$ for that $k$. 
Since all of the $P_k$ are conjugate in $\cC_n$ by Remark~\ref{swap}, we have $\chi(P_k)=-1$ {\em for all} $k$; for the same reason we have $\chi(H_k)=-1$ for all $k$. The defining relation $(R10)$ forces $\chi(Z_j)=-1$ for all $j$, but then $(B2)$ would give a contradiction. Therefore $\chi(P_k)=1$ for all $k$. A similar argument shows that $\chi(H_k)=1$ for all $k$. Then $(B2)$ forces $\chi(Z_j)=1$ for all $j$ and hence $\chi$ must be the trivial character.
\end{proof}

\begin{remark}
Theorem \ref{c=cc} can also be deduced from the fact that $Sp(2n,2)$ is a perfect group except for $n=1,2$ \cite{BRW,He}. It can also be deduced from Theorem \ref{normalc}, see Remark~\ref{altpf}. Our proof above does not require any knowledge from $Sp(2n,2)$.
\end{remark}

Some patterns can be observed from our previous results, suggesting the following results to hold for general $n$.
In fact, they are consequences of known results about $Sp(2n,2)$ \cite{BRW}, 
but we believe that they should be explicitly referred to properties of $\cC_n$. 
Recall that $\cP_n$ is normal in $\cC_n$. In other words, $\cC_n$ acts on $\cP_n$ by conjugation, which is in fact a transitive action. 
\begin{lemma}\cite{GAE,Zh}\label{trac}
$\cC_n$ acts on $\cP_n$ transitively by conjugation. Hence $\cP_n$ does not contain any non-trivial proper subgroup which is normal in $\cC_n$.
\end{lemma}

\begin{theorem}\label{normalc}
When $n\geq 3$, the $n$-qubit Pauli group $\cP_n$ is the only non-trivial proper normal subgroup of $\cC_n$.  
\end{theorem}
\begin{proof}
Let $N$ be a non-trivial proper normal subgroup of $\cC_n$ with $N\neq \cP_n$. 
Denote by $\pi:\cC_n\rightarrow \cC_n/\cP_n$ the natural quotient map, then $\pi(N)$ is a normal subgroup of $\cC_n/\cP_n$.
Since $\cC_n/\cP_n\cong Sp(2n,2)$ is simple, it is impossible to have $\cP_n \subsetneq N$, or $\pi(N)$ would be a non-trivial proper normal subgroup in $\cC_n/\cP_n$.
Thus we must have $\cP_n \nsubseteq N$. 

Suppose $\pi(N)=\cP_n$, which implies that $N\subseteq \cP_n$. Then $N$ is a non-trivial proper subgroup of $\cP_n$ which is normal in $\cC_n$, but this contradicts to Lemma \ref{trac}.

Suppose now $\pi(N)=\cC_n/\cP_n$. Then $\cC_n=N\cP_n$ with $\cP_n$, $N$ both normal in $\cC_n$. The intersection $N\cap \cP_n$ is a subgroup of $\cP_n$ which is normal in $\cC_n$, and hence $N\cap \cP_n$ must be the trivial group by Lemma \ref{trac}. Thus $\cC_n=N\times \cP_n$, and in particular we have $N$ and $\cP_n$ commute with each other in $\cC_n$. But then the conjugation of $\cC_n$ on $\cP_n$ would be the trivial action since $\cP_n$ is an abelian group, contradicts to Lemma \ref{trac}.
\end{proof}

\begin{remark}\label{altpf}
As a consequence of Theorem \ref{normalc}, for $n\geq 3$, the only quotient groups of $\cC_n$ are the trivial group, $\cC_n$, and $\cC_n/\cP_n\cong Sp(2n,2)$. This gives an alternative proof of Theorem \ref{c=cc}. 
\end{remark}

\begin{corollary}\label{faith}
For $n\geq 3$, any representation of $\cC_n$ which is non-trivial on $\cP_n$ is faithful.
In particular, the matrix representation $\cM_{2^n}$ for $\cC_n$ is faithful for all $n$.
\end{corollary}
\begin{proof}
The first statement is a consequence of Theorem \ref{normalc}, which implies the second statement for $n\geq 3$. The cases $n=1,2$ can be observed from our previous results.
\end{proof}

\begin{remark}\label{conj}
Theorem \ref{normalc} and Corollary \ref{faith} were only conjectures in a previous draft of this article.
We thank Jason Saied for providing his comment \cite{Sa} that turn them into affirmative results.
\end{remark}

\begin{remark}\label{cntab}
In a very recent article, \cite{Ma} developed the Clifford theory for $\cC_n$. In particular, applying \cite[Theorem 5.3]{Ma}, one can lift an irreducible character of $\cC_{n-1}$ to obtain an irreducible one of $\cC_{n}$. Hence one can obtain a part of irreducible characters of $\cC_{n}$ from the character table of $\cC_{n-1}$. According to \cite[Theorem 4.4]{Ma}, the remaining irreducible characters of $\cC_n$ are precisely those lifted from the irreducible character of $Sp(2n,2)$. By this fashion, one should be able to obtain the complete character table of $\cC_n$ for all $n$, assuming the character table of $Sp(2n,2)$ is available.
\end{remark}

\begin{remark}\label{intab}
In a previous draft of this article we proposed the following conjecture: the character table of $\cC_n$ consists only integers. However, this is true only for $1\leq n\leq 3$. When $n=4$, it is a fact that there exist irreducible characters of $Sp(8,2)$ containing non-integral values. According to \cite[Lemma 4.6, Lemma 5.1]{Ma}, these irreducible characters of $Sp(8,2)$ lift to irreducible characters of $\cC_4$ with non-integral values and they keep to exist for all $\cC_n$ with $n\geq 5$.
\end{remark}

To conclude this article, we highlight an interesting similarity between the symmetric groups and the Clifford groups, which we don't know is it a coincidence or not.
It is well known that the representation theory of the symmetric groups can be constructed combinatorially by Young tableaux \cite{FH, Sag}. 
It will be great if one can develop a parallel combinatorial object to study the representation theory of $\cC_n$. 
\begin{proposition}
The following facts for $S_n$ hold for all $n\geq 2$:
\begin{enumerate}
\item The standard representation $\mathbb{C}^n$ decomposes multiplicity-freely as the direct sum of the trivial representation and an irreducible representation of codimension 1.
\item The standard representation $\mathbb{C}^n$ is a faithful representation.
\item \cite[Theorem~19.10]{JL} Any irreducible representation of $S_n$ must appear as an irreducible constituent of $(\mathbb{C}^n)^{\otimes d}$ for some $d$ large enough. 
\end{enumerate}
\end{proposition}

\begin{proposition}
The following facts for $\cC_n$ hold for all $n\geq 1$:
\begin{enumerate}
\item \cite{WF} The matrix representation $\cM_{2^n}$ decomposes multiplicity-freely as the direct sum of the trivial representation and an irreducible representation of codimension 1.
\item (Corollary \ref{faith}) The matrix representation $\cM_{2^n}$ is a faithful representation.
\item \cite[Theorem~19.10]{JL} Any irreducible representation of $\cC_n$ must appear as an irreducible constituent of $(\cM_{2^n})^{\otimes d}$ for some $d$ large enough.
\end{enumerate}
\end{proposition}

\section{Appendix}
\subsection{Some GAP codes}
We provide some codes for the computation performed by GAP as follows.
Note that although Theorem \ref{presc} provides an explicit way to define the group, it is not easy for GAP to calculate the desired results even for $n=2$. To overcome this obstacle, we embed our group into a larger symmetric group on which the calculation can be done much faster than using the abstract symbols and defining relations.

\begin{verbatim}
gap>f:=FreeGroup("h1","h2","p1","p2","z");;
gap>AssignGeneratorVariables(f);;
gap>rels:=[ z^2,h1^2,h2^2,p1^4,p2^4,
(h1*p1)^3,(h2*p2)^3,
(h1*p1*h1*p1^3)^3,(h2*p2*h2*p2^3)^3,
(h1*p1^2*h1*p1^2)^2,(h2*p2^2*h2*p2^2)^2,
Comm(h1,h2),Comm(p1,p2),
Comm(z,p1),Comm(z,p2),
z*h1*p1^2*h1/(h1*p1^2*p2^2*h1*z),
z*h2*p2^2*h2/(h2*p1^2*p2^2*h2*z),
z*h1*z/(p1*h1*p1*p2*z*h1*p1),
z*h2*z/(p2*h2*p1*p2*z*h2*p2)];;
gap>g:=f/rels;;
gap>AssignGeneratorVariables(g);;
gap>IsomorphismPermGroup(g);
\end{verbatim}

\subsection{A presentation of $Sp(2n,2)$}
As a consequence of Theorem \ref{presc} and the fact that $\cC_n/\cP_n\cong Sp(2n,2)$, one can write down a presentation for the symplectic group $Sp(2n,2)$ as follows.
\begin{theorem}\label{presp}
The symplectic group $Sp(2n,2)$ is generated by the symbols 
$$\{H_i, P_i, Z_{j} \,\, | \,\, 1\leq i\leq n, \, 1\leq j\leq n-1\}$$ subjected to all relations in Theorem \ref{presc} together with following extra relations:
\begin{itemize}
\item(SP1) $P_i^2=1$ for all $1\leq i\leq n$. \\[-3mm]
\item(SP2) $H_iP_i^3(H_iP_i)^2=1$ for all $1\leq i\leq n$.  \\[-3mm]
\item(SP3) $(H_iP_i)^2P_i^2H_iP_i=1$ for all $1\leq i\leq n$.  \\[-3mm]
\end{itemize}
\end{theorem}
\begin{proof}
Since $Sp(2n,2)\cong\cC_n/\cP_n$, we only need to impose extra relations sending precisely the $n$-qubit Pauli group to 1. Direct computation shows that $P^2=Z$, $HP^3 (HP)^2=X$ and $(HP)^2 P^2 HP=Y$ in $\cC_1$, and hence the subgroup generated by $P_i^2$, $H_iP_i^3(H_iP_i)^2$ and $(H_iP_i)^2P_i^2H_iP_i$ for all $1\leq i\leq n$ is exactly the $n$-qubit Pauli group $\cP_n$.
\end{proof}
\begin{remark}
Some relations are in fact redundant. 
For example, it is enough to set only $i=1$ in $(SP1)$--$(SP3)$, where the swapping operator in Remark \ref{swap} will bring the relations for all other $i$. Also, after applying $(SP1)$, both of $(SP2)$ and $(SP3)$ can be further reduced to the single relation $(H_iP_i)^3=1$.
\end{remark}

\subsection{Character table of Three-qubit Clifford group}\label{appc3}
Applying Theorem \ref{presc}, the following character table of $\cC_3$ is calculated by GAP.\\[5mm]

\setlength\tabcolsep{1pt}
\footnotesize
\begin{tabular}{*{34}{|c}}
\hline
 & 1 & 2&3&4&5&6&7&8&9&10&11&12&13&14&15&16&17&18&19&20&21&22&23&24&25&26&27&28&29&30&31&32\\
\hline\hline
$\chi_{1}$&1& 1& 1& 1& 1& 1& 1& 1& 1& 1& 1& 1& 1& 1& 1& 1& 1& 1& 1& 1& 1& 1& 1& 1& 1& 1& 1& 1& 1& 1& 1& 1\\ \hline
$\chi_{2}$&7& 7& -2& 1& 3& 3& 3& 3& 4& 4& -3& -3& -3& 0& 0& 0& 0& 0& 1& 1& 1& -1& -1& -1& 2& 2& 1& 1& -1& -1& -5& -5\\ \hline
$\chi_{3}$&15& 15& -3& 0& 3& 3& 3& 3& 0& 0& 1& 1& 1& 0& 0& 0& -2& -2& -3& -3& -3& 7& 7& 7& -2& -2& 3& 3& 1& 1& -5& -5\\ \hline
$\chi_{4}$&21& 21& 3& 0& 5& 5& 5& 5& 6& 6& -3& -3& -3& 2& 2& 2& 0& 0& -3& -3& -3& 5& 5& 5& 2& 2& 0& 0& 2& 2& -11& -11\\ \hline
$\chi_{5}$&21& 21& 3& 0& 1& 1& 1& 1& 6& 6& 3& 3& 3& -2& -2& -2& 0& 0& -1& -1& -1& -3& -3& -3& 0& 0& 0& 0& 0& 0& 9& 9\\ \hline
$\chi_{6}$&27& 27& 0& 0& 7& 7& 7& 7& 9& 9& 5& 5& 5& 1& 1& 1& -1& -1& 1& 1& 1& 3& 3& 3& 3& 3& 0& 0& 0& 0& 15& 15\\ \hline
$\chi_{7}$&35& 35& -1& -1& 7& 7& 7& 7& 5& 5& 1& 1& 1& 1& 1& 1& 1& 1& 5& 5& 5& 11& 11& 11& -1& -1& 2& 2& 2& 2& 15& 15\\ \hline
$\chi_{8}$&35& 35& -1& -1& -5& -5& -5& -5& 5& 5& -1& -1& -1& 1& 1& 1& -1& -1& -1& -1& -1& 3& 3& 3& -3& -3& 2& 2& 0& 0& -5& -5\\ \hline
$\chi_{9}$&56& 56& 2& -1& 8& 8& 8& 8& 11& 11& -4& -4& -4& -1& -1& -1& -1& -1& 4& 4& 4& -8& -8& -8& 1& 1& 2& 2& -2& -2& -24& -24\\ \hline
$\chi_{10}$&63& -1& 0& 0& -1& -1& -1& 15& -1& 15& -1& -1& 7& 3& -1& -1& -1& 1& -1& -1& 7& -1& 15& -1& -1& 3& 3& -1& -1& 3& -1& 31\\ \hline
$\chi_{11}$&63& -1& 0& 0& -1& -5& 3& 11& -1& 15& 1& 1& -7& -1& -1& 3& 1& -1& 1& -3& 5& -1& -9& 7& -1& 3& 3& -1& 1& -3& -1& -29\\ \hline
$\chi_{12}$&70& 70& 7& 1& 6& 6& 6& 6& -5& -5& 2& 2& 2& 3& 3& 3& -1& -1& 2& 2& 2& -10& -10& -10& -1& -1& 1& 1& -1& -1& -10& -10\\ \hline
$\chi_{13}$&84& 84& 3& 0& 4& 4& 4& 4& -6& -6& 0& 0& 0& -2& -2& -2& 0& 0& 0& 0& 0& 20& 20& 20& 2& 2& 3& 3& -1& -1& 4& 4\\ \hline
$\chi_{14}$&105& 105& -3& 0& 5& 5& 5& 5& 15& 15& -5& -5& -5& -1& -1& -1& 1& 1& -1& -1& -1& 1& 1& 1& 1& 1& -3& -3& 1& 1& -35& -35\\ \hline
$\chi_{15}$&105& 105& 6& 0& 9& 9& 9& 9& 0& 0& -3& -3& -3& 0& 0& 0& 0& 0& -3& -3& -3& -7& -7& -7& -4& -4& 3& 3& -1& -1& 25& 25\\
\hline
 $\chi_{16}$&105& 105& 6& 0& -3& -3& -3& -3& 0& 0& -1& -1& -1& 0& 0& 0& 2& 2& 3& 3& 3& 17& 17& 17& 2& 2& 3& 3& -1& -1& 5& 5\\ \hline
 $\chi_{17}$&120& 120& -6& 0& 8& 8& 8& 8& 15& 15& 4& 4& 4& -1& -1& -1& 1& 1& -4& -4& -4& -8& -8& -8& 1& 1& 0& 0& -2& -2& 40& 40\\ \hline
\end{tabular}\\[.2cm]

\setlength\tabcolsep{1pt}
\footnotesize
\begin{tabular}{*{40}{|c}}
\hline
&33&34&35&36&37&38&39&40&41&42&43&44&45&46&47&48&49&50&51&52&53&54&55&56&57&58&59&60&61&62&63&64&65&66&67 \\
\hline\hline
$\chi_{1}$&1& 1& 1& 1& 1& 1& 1& 1& 1& 1& 1& 1& 1& 1& 1& 1& 1& 1& 1& 1& 1& 1& 1& 1& 1& 1& 1& 1& 1& 1& 1& 1& 1& 1& 1\\ \hline
$\chi_{2}$&-5& 1& 1& -1& -1& -1& -1& -1& -1& 3& 3& 2& 0& 2& 2& 0& 0& -1& -1& -1& -1& -1& 1& 1& 1& 1& 1& 1& -1& -2& -2& -2& -2& -2& 0\\ \hline
$\chi_{3}$&-5& 1& 1& -1& -1& -1& -1& -1& -1& -1& -1& 1& -1& 0& 0& 0& 0& 1& 1& 3& 3& 3& -1& -1& 1& 1& 1& 1& 0& -2& -2& -2& 0& 0& 1\\ \hline
$\chi_{4}$&-11& -2& -2& -3& -3& -3& -3& 0& 0& 1& 1& -1& 1& 1& 1& -1& -1& -1& -1& 1& 1& 1& -1& -1& 1& 1& 1& 1& 1& -2& -2& -2& 0& 0& 0\\ \hline
$\chi_{5}$&9& 0& 0& -3& -3& -3& -3& 0& 0& 5& 5& 3& -1& 1& 1& -1& -1& -1& -1& 1& 1& 1& 1& 1& -1& -1& -1& -1& 1& 0& 0& 0& 2& 2& 0\\ \hline
$\chi_{6}$&15& 0& 0& 3& 3& 3& 3& 0& 0& 3& 3& 0& 0& 2& 2& 0& 0& 1& 1& -1& -1& -1& -1& -1& 1& 1& 1& 1& -1& 3& 3& 3& 1& 1& -1\\ \hline
$\chi_{7}$&15& 0& 0& 3& 3& 3& 3& 0& 0& -1& -1& -1& -1& 0& 0& 0& 0& -1& -1& 3& 3& 3& 1& 1& 1& 1& 1& 1& 0& 3& 3& 3& -1& -1& 0\\ \hline
$\chi_{8}$&-5& -2& -2& 3& 3& 3& 3& 0& 0& 7& 7& 3& 1& 0& 0& 0& 0& 1& 1& -1& -1& -1& 1& 1& -1& -1& -1& -1& 0& 1& 1& 1& -1& -1& 0\\ \hline
$\chi_{9}$&-24& 0& 0& 0& 0& 0& 0& 0& 0& 0& 0& -2& 0& 1& 1& 1& 1& 0& 0& 0& 0& 0& 0& 0& 0& 0& 0& 0& 1& -3& -3& -3& 1& 1& 0\\ \hline
$\chi_{10}$&-1& -1& 1& -1& 7& -1& -1& -1& 1& 3& -1& 0& 0& -1& 3& 1& -1& 1& -1& 3& -1& -1& -1& 1& -1& 3& -1& -1& 0& -1& -1& 7& 1& -1& 0\\ \hline
$\chi_{11}$&3& -1& 1& -1& -1& 3& -5& 1& -1& 3& -1& 0& 0& -1& 3& 1& -1& -1& 1& -1& -1& 3& -1& 1& 1& 1& 1& -3& 0& -1& 3& -5& -1& 1& 0\\ \hline
$\chi_{12}$&-10& -1& -1& -2& -2& -2& -2& 1& 1& 2& 2& -1& -1& 0& 0& 0& 0& 0& 0& 2& 2& 2& 0& 0& -2& -2& -2& -2& 0& -1& -1& -1& -1& -1& 0\\ \hline
$\chi_{13}$&4& 1& 1& 4& 4& 4& 4& 1& 1& 4& 4& -1& 1& -1& -1& -1& -1& 0& 0& 4& 4& 4& 0& 0& 0& 0& 0& 0& -1& -2& -2& -2& 0& 0& 0\\ \hline
$\chi_{14}$&-35& 1& 1& 1& 1& 1& 1& 1& 1& 5& 5& 1& -1& 0& 0& 0& 0& 1& 1& 1& 1& 1& -1& -1& -1& -1& -1& -1& 0& 1& 1& 1& -1& -1& 0\\ \hline
$\chi_{15}$&25& 1& 1& 1& 1& 1& 1& 1& 1& -3& -3& 2& 0& 0& 0& 0& 0& -1& -1& -3& -3& -3& -1& -1& 1& 1& 1& 1& 0& 4& 4& 4& 0& 0& 0\\ \hline
$\chi_{16}$&5& -1& -1& -7& -7& -7& -7& -1& -1& -3& -3& 2& 0& 0& 0& 0& 0& 1& 1& 1& 1& 1& -1& -1& -1& -1& -1& -1& 0& 2& 2& 2& 0& 0& 0\\ \hline
$\chi_{17}$&40& -2& -2& 0& 0& 0& 0& 0& 0& 0& 0& -2& 0& 0& 0& 0& 0& 0& 0& 0& 0& 0& 0& 0& 0& 0& 0& 0& 0& 1& 1& 1& -1& -1& 1\\ \hline
\end{tabular}\\[.2cm]

\setlength\tabcolsep{1pt}
\footnotesize
\begin{tabular}{*{34}{|c}}
\hline
 & 1 & 2&3&4&5&6&7&8&9&10&11&12&13&14&15&16&17&18&19&20&21&22&23&24&25&26&27&28&29&30&31&32\\
 \hline\hline
 $\chi_{18}$&168& 168& 6& 0& 8& 8& 8& 8& 6& 6& 0& 0& 0& 2& 2& 2& 0& 0& 0& 0& 0& 8& 8& 8& 2& 2& -3& -3& -1& -1& 40& 40\\ \hline
$\chi_{19}$&189& 189& 0& 0& 1& 1& 1& 1& 9& 9& -1& -1& -1& 1& 1& 1& -1& -1& -5& -5& -5& 21& 21& 21& 3& 3& 0& 0& 0& 0& -39& -39\\ \hline
$\chi_{20}$&189& 189& 0& 0& 13& 13& 13& 13& 9& 9& 1& 1& 1& 1& 1& 1& 1& 1& 1& 1& 1& -3& -3& -3& -3& -3& 0& 0& 0& 0& -51& -51\\ \hline
$\chi_{21}$&189& 189& 0& 0& -11& -11& -11& -11& 9& 9& 1& 1& 1& 1& 1& 1& 1& 1& 1& 1& 1& -3& -3& -3& -3& -3& 0& 0& 0& 0& 21& 21\\ \hline
$\chi_{22}$&210& 210& 3& 0& 2& 2& 2& 2& 15& 15& 2& 2& 2& -1& -1& -1& -1& -1& 2& 2& 2& 2& 2& 2& -1& -1& 0& 0& 2& 2& 50& 50\\ \hline
$\chi_{23}$&210& 210& -6& 0& 10& 10& 10& 10& -15& -15& -2& -2& -2& 1& 1& 1& 1& 1& -2& -2& -2& -14& -14& -14& 1& 1& 3& 3& 1& 1& 10& 10\\ \hline
$\chi_{24}$&216& 216& 0& 0& 8& 8& 8& 8& -9& -9& 4& 4& 4& -1& -1& -1& 1& 1& -4& -4& -4& 24& 24& 24& -3& -3& 0& 0& 0& 0& -24& -24\\ \hline
$\chi_{25}$&280& 280& -8& 1& 8& 8& 8& 8& -5& -5& -4& -4& -4& -1& -1& -1& -1& -1& 4& 4& 4& 24& 24& 24& -3& -3& -2& -2& 0& 0& 40& 40\\ \hline
$\chi_{26}$&280& 280& 10& 1& -8& -8& -8& -8& 10& 10& 0& 0& 0& -2& -2& -2& 0& 0& 0& 0& 0& -8& -8& -8& -2& -2& 1& 1& 1& 1& -40& -40\\ \hline
 $\chi_{27}$&315& -5& 0& 0& -5& -5& 11& 11& -2& 30& 3& -1& -9& 2& -2& 2& 0& 0& 3& -1& -9& -5& 27& 11& -2& 6& -3& 1& -1& 3& -5& -85\\ \hline
$\chi_{28}$&315& -5& 0& 0& -5& 7& -1& 23& -2& 30& -3& 1& 9& 2& -2& 2& 0& 0& -3& 5& -3& -5& 3& 19& -2& 6& -3& 1& 1& -3& -5& 95\\ \hline
 $\chi_{29}$&315& 315& -9& 0& 3& 3& 3& 3& 0& 0& 3& 3& 3& 0& 0& 0& 0& 0& 3& 3& 3& -21& -21& -21& 0& 0& 0& 0& 0& 0& -45& -45\\ \hline
$\chi_{30}$&315& -5& 0& 0& -5& -1& 7& 15& 1& -15& 1& -3& 5& -3& 1& 1& 1& -1& 1& 1& -7& -5& 51& 3& 1& -3& 6& -2& 0& 0& -5& -25\\ \hline
 $\chi_{31}$&315& -5& 0& 0& -5& 3& 3& 19& 1& -15& -1& 3& -5& 1& 1& -3& -1& 1& -1& 3& -5& -5& -21& 27& 1& -3& 6& -2& 0& 0& -5& 35\\ \hline
$\chi_{32}$&336& 336& -6& 0& -16& -16& -16& -16& 6& 6& 0& 0& 0& 2& 2& 2& 0& 0& 0& 0& 0& 16& 16& 16& -2& -2& 0& 0& -2& -2& -16& -16\\ \hline
$\chi_{33}$&378& 378& 0& 0& 2& 2& 2& 2& -9& -9& 2& 2& 2& -1& -1& -1& -1& -1& 2& 2& 2& -6& -6& -6& 3& 3& 0& 0& 0& 0& -30& -30\\ \hline
$\chi_{34}$&378& -6& 0& 0& -2& 2& -6& 34& -3& 45& 2& 2& -14& 1& 1& -3& -1& 1& -2& 2& 2& 2& -6& -6& -1& 3& 0& 0& 0& 0& 6& -126\\ \hline
$\chi_{35}$&378& -6& 0& 0& -2& -14& 10& 18& -3& 45& -2& -2& 14& -3& 1& 1& 1& -1& 2& -2& -2& 2& -6& -6& -1& 3& 0& 0& 0& 0& 6& 114\\ \hline
$\chi_{36}$&405& 405& 0& 0& -3& -3& -3& -3& 0& 0& -3& -3& -3& 0& 0& 0& 0& 0& -3& -3& -3& -27& -27& -27& 0& 0& 0& 0& 0& 0& 45& 45\\ \hline
$\chi_{37}$&420& 420& -3& 0& -12& -12& -12& -12& 0& 0& 0& 0& 0& 0& 0& 0& 0& 0& 0& 0& 0& 4& 4& 4& 4& 4& 3& 3& 1& 1& 20& 20\\ \hline
$\chi_{38}$&512& 512& 8& -1& 0& 0& 0& 0& -16& -16& 0& 0& 0& 0& 0& 0& 0& 0& 0& 0& 0& 0& 0& 0& 0& 0& -4& -4& 0& 0& 0& 0\\ \hline
$\chi_{39}$&567& -9& 0& 0& -9& -1& 15& 15& 0& 0& 3& -5& 3& 0& 0& 0& 0& 0& 3& -5& 3& -9& -9& 39& 0& 0& 0& 0& 0& 0& -9& -81\\ \hline
$\chi_{40}$&567& -9& 0& 0& -9& 11& 3& 27& 0& 0& -3& 5& -3& 0& 0& 0& 0& 0& -3& 1& 9& -9& 63& 15& 0& 0& 0& 0& 0& 0& -9& 99\\ \hline
$\chi_{41}$&630& -10& 0& 0& 2& -2& -10& 30& -1& 15& -2& 2& 2& 3& -1& -1& 1& -1& 2& 2& -14& -2& 54& -10& 1& -3& 3& -1& -1& 3& 10& -130\\ \hline
$\chi_{42}$&630& -10& 0& 0& 2& -2& -10& 30& -1& 15& -2& 2& 2& 3& -1& -1& 1& -1& 2& -6& 10& -2& -42& 22& 1& -3& 3& -1& 1& -3& 10& -130\\ \hline
\end{tabular}\\[.2cm]

\setlength\tabcolsep{1pt}
\footnotesize
\begin{tabular}{*{40}{|c}}
\hline
&33&34&35&36&37&38&39&40&41&42&43&44&45&46&47&48&49&50&51&52&53&54&55&56&57&58&59&60&61&62&63&64&65&66&67 \\
\hline\hline
$\chi_{18}$&40& 1& 1& 8& 8& 8& 8& -1& -1& 0& 0& 2& 0& -2& -2& 0& 0& 0& 0& 0& 0& 0& 0& 0& 0& 0& 0& 0& 1& -2& -2& -2& 0& 0& 0\\ \hline
$\chi_{19}$&-39& 0& 0& -3& -3& -3& -3& 0& 0& -3& -3& 0& 0& -1& -1& 1& 1& -1& -1& 1& 1& 1& 1& 1& -1& -1& -1& -1& -1& 3& 3& 3& 1& 1& 0\\ \hline
$\chi_{20}$&-51& 0& 0& -3& -3& -3& -3& 0& 0& -3& -3& 0& 0& -1& -1& -1& -1& 1& 1& -3& -3& -3& 1& 1& 1& 1& 1& 1& -1& -3& -3& -3& 1& 1& 0\\ \hline
$\chi_{21}$&21& 0& 0& -3& -3& -3& -3& 0& 0& 9& 9& 0& 0& -1& -1& 1& 1& -1& -1& 1& 1& 1& -1& -1& 1& 1& 1& 1& -1& -3& -3& -3& 1& 1& 0\\ \hline
$\chi_{22}$&50& 2& 2& -6& -6& -6& -6& 0& 0& -2& -2& -1& 1& 0& 0& 0& 0& 0& 0& -2& -2& -2& 0& 0& -2& -2& -2& -2& 0& -1& -1& -1& -1& -1& 0\\ \hline
$\chi_{23}$&10& 1& 1& 2& 2& 2& 2& -1& -1& 6& 6& -2& 0& 0& 0& 0& 0& 0& 0& -2& -2& -2& 0& 0& -2& -2& -2& -2& 0& 1& 1& 1& 1& 1& 0\\ \hline
$\chi_{24}$&-24& 0& 0& 0& 0& 0& 0& 0& 0& 0& 0& 0& 0& 1& 1& 1& 1& 0& 0& 0& 0& 0& 0& 0& 0& 0& 0& 0& 1& -3& -3& -3& -1& -1& -1\\ \hline
$\chi_{25}$&40& -2& -2& 0& 0& 0& 0& 0& 0& 0& 0& 0& 0& 0& 0& 0& 0& 0& 0& 0& 0& 0& 0& 0& 0& 0& 0& 0& 0& 1& 1& 1& 1& 1& 0\\ \hline
$\chi_{26}$&-40& -1& -1& 8& 8& 8& 8& -1& -1& 0& 0& -2& 0& 0& 0& 0& 0& 0& 0& 0& 0& 0& 0& 0& 0& 0& 0& 0& 0& 2& 2& 2& 0& 0& 0\\ \hline
$\chi_{27}$&11& 1& -1& 3& -5& 3& -13& -1& 1& 3& -1& 0& 0& 0& 0& 0& 0& -1& 1& 3& -1& -1& 1& -1& 3& -1& -1& -1& 0& -2& 2& 2& 0& 0& 0\\ \hline
$\chi_{28}$&-1& 1& -1& -5& 11& -1& 7& 1& -1& 3& -1& 0& 0& 0& 0& 0& 0& 1& -1& -1& -1& 3& 1& -1& 1& 1& -3& 1& 0& -2& 2& 2& 0& 0& 0\\ \hline
$\chi_{29}$&-45& 0& 0& 3& 3& 3& 3& 0& 0& -5& -5& 3& 1& 0& 0& 0& 0& -1& -1& 3& 3& 3& -1& -1& -1& -1& -1& -1& 0& 0& 0& 0& 0& 0& 0\\ \hline
$\chi_{30}$&7& -2& 2& 3& 3& -1& -9& 0& 0& 3& -1& 0& 0& 0& 0& 0& 0& 1& -1& 7& -1& -5& 1& -1& 1& 1& -3& 1& 0& 1& 1& -7& -1& 1& 0\\ \hline
$\chi_{31}$&3& -2& 2& -5& 3& 3& 3& 0& 0& 3& -1& 0& 0& 0& 0& 0& 0& -1& 1& -5& -1& 7& 1& -1& 3& -1& -1& -1& 0& 1& -3& 5& 1& -1& 0\\ \hline
$\chi_{32}$&-16& 2& 2& 0& 0& 0& 0& 0& 0& 0& 0& -2& 0& 1& 1& -1& -1& 0& 0& 0& 0& 0& 0& 0& 0& 0& 0& 0& 1& 2& 2& 2& 0& 0& 0\\ \hline
$\chi_{33}$&-30& 0& 0& -6& -6& -6& -6& 0& 0& 6& 6& 0& 0& -2& -2& 0& 0& 0& 0& -2& -2& -2& 0& 0& 2& 2& 2& 2& 1& 3& 3& 3& -1& -1& 0\\ \hline
$\chi_{34}$&2& 0& 0& 2& -6& -2& 6& 0& 0& 6& -2& 0& 0& -1& 3& -1& 1& 0& 0& -2& 2& -2& 0& 0& -2& 2& -2& 2& 0& 3& -1& -9& -1& 1& 0\\ \hline
$\chi_{35}$&-14& 0& 0& 2& -6& -2& 6& 0& 0& 6& -2& 0& 0& -1& 3& -1& 1& 0& 0& -2& 2& -2& 0& 0& 2& -2& 2& -2& 0& 3& -5& 3& 1& -1& 0\\ \hline
$\chi_{36}$&45& 0& 0& -3& -3& -3& -3& 0& 0& -3& -3& 0& 0& 0& 0& 0& 0& 1& 1& 5& 5& 5& 1& 1& 1& 1& 1& 1& 0& 0& 0& 0& 0& 0& -1\\ \hline
$\chi_{37}$&20& -1& -1& 4& 4& 4& 4& 1& 1& -4& -4& 1& -1& 0& 0& 0& 0& 0& 0& -4& -4& -4& 0& 0& 0& 0& 0& 0& 0& -4& -4& -4& 0& 0& 0\\ \hline
$\chi_{38}$&0& 0& 0& 0& 0& 0& 0& 0& 0& 0& 0& 0& 0& 2& 2& 0& 0& 0& 0& 0& 0& 0& 0& 0& 0& 0& 0& 0& -1& 0& 0& 0& 0& 0& 1\\ \hline
$\chi_{39}$&15& 0& 0& -1& -9& 7& -9& 0& 0& 3& -1& 0& 0& 1& -3& -1& 1& 1& -1& -5& -1& 7& -1& 1& -1& 3& -1& -1& 0& 0& 0& 0& 0& 0& 0\\ \hline
$\chi_{40}$&3& 0& 0& -1& 15& -5& 3& 0& 0& 3& -1& 0& 0& 1& -3& -1& 1& -1& 1& 7& -1& -5& -1& 1& 1& 1& 1& -3& 0& 0& 0& 0& 0& 0& 0\\ \hline
$\chi_{41}$&-2& 1& -1& -2& -10& 2& 10& 1& -1& -6& 2& 0& 0& 0& 0& 0& 0& 0& 0& 2& -2& 2& 0& 0& -2& 2& 2& -2& 0& 1& 1& -7& 1& -1& 0\\ \hline
$\chi_{42}$&-2& 1& -1& 6& -2& -6& 2& -1& 1& -6& 2& 0& 0& 0& 0& 0& 0& 0& 0& 2& -2& 2& 0& 0& 2& -2& -2& 2& 0& 1& 1& -7& 1& -1& 0\\ \hline
\end{tabular}\\[.2cm]

\setlength\tabcolsep{1pt}
\footnotesize
\begin{tabular}{*{34}{|c}}
\hline
 & 1 & 2&3&4&5&6&7&8&9&10&11&12&13&14&15&16&17&18&19&20&21&22&23&24&25&26&27&28&29&30&31&32\\
\hline\hline
$\chi_{43}$&630& -10& 0& 0& -10& 14& 14& -18& -1& 15& 2& -2& -2& 3& -1& -1& -1& 1& 2& -2& -2& -10& 54& 22& -1& 3& 3& -1& 1& -3& -10& -50\\ \hline
$\chi_{44}$&630& -10& 0& 0& 2& -18& 6& 14& -1& 15& 2& -2& -2& -1& -1& 3& -1& 1& -2& -2& 14& -2& 54& -10& 1& -3& 3& -1& -1& 3& 10& 110\\ \hline
$\chi_{45}$&630& -10& 0& 0& 2& -18& 6& 14& -1& 15& 2& -2& -2& -1& -1& 3& -1& 1& -2& 6& -10& -2& -42& 22& 1& -3& 3& -1& 1& -3& 10& 110\\ \hline
$\chi_{46}$&630& -10& 0& 0& -10& 22& 6& -10& -1& 15& -2& 2& 2& -1& -1& 3& 1& -1& -2& 2& 2& -10& 6& 38& -1& 3& 3& -1& -1& 3& -10& 70\\ \hline
$\chi_{47}$&945& -15& 0& 0& 1& 1& -15& 49& -3& 45& 1& -3& 5& 1& 1& -3& 1& -1& 1& -3& 5& 1& 33& -15& 1& -3& 0& 0& 0& 0& -15& 225\\ \hline
$\chi_{48}$&945& -15& 0& 0& 1& -7& 9& -23& -3& 45& 1& 1& -7& 1& 1& -3& 1& -1& 1& 1& -7& 1& 33& -15& 1& -3& 0& 0& 0& 0& -15& -135\\ \hline
$\chi_{49}$&945& -15& 0& 0& 1& -27& 13& 21& -3& 45& -1& 3& -5& -3& 1& 1& -1& 1& -1& -1& 7& 1& -39& 9& 1& -3& 0& 0& 0& 0& -15& -195\\ \hline
$\chi_{50}$&945& -15& 0& 0& 1& 13& -11& -3& -3& 45& -1& -1& 7& -3& 1& 1& -1& 1& -1& 3& -5& 1& -39& 9& 1& -3& 0& 0& 0& 0& -15& 165\\ \hline
$\chi_{51}$&1008& -16& 0& 0& -16& 16& 16& 16& 2& -30& 0& 0& 0& -2& 2& -2& 0& 0& 0& 0& 0& -16& 48& 48& 2& -6& -6& 2& 0& 0& -16& 16\\ \hline
$\chi_{52}$&1260& -20& 0& 0& 4& 12& -4& -52& -2& 30& 0& 0& 0& 2& -2& 2& 0& 0& 0& 0& 0& -4& 12& 12& 2& -6& 6& -2& 0& 0& 20& -20\\ \hline
$\chi_{53}$&1512& -24& 0& 0& -8& -8& 24& -8& -3& 45& -4& 4& 4& 1& 1& -3& -1& 1& 4& -4& -4& 8& -24& -24& -1& 3& 0& 0& 0& 0& 24& 216\\ \hline
$\chi_{54}$&1512& -24& 0& 0& -8& 24& -8& 24& -3& 45& 4& -4& -4& -3& 1& 1& 1& -1& -4& 4& 4& 8& -24& -24& -1& 3& 0& 0& 0& 0& 24& -264\\ \hline
$\chi_{55}$&1890& -30& 0& 0& -10& 10& 2& 42& 3& -45& 2& -6& 10& 3& -1& -1& -1& 1& -2& 2& 2& 10& -30& -30& 1& -3& 0& 0& 0& 0& 30& -150\\ \hline
$\chi_{56}$&1890& -30& 0& 0& 2& -14& 2& 18& 3& -45& -2& 2& 2& 3& -1& -1& 1& -1& -2& 2& 2& 2& -78& 18& -1& 3& 0& 0& 0& 0& -30& -30\\ \hline
$\chi_{57}$&1890& -30& 0& 0& -10& -6& 18& 26& 3& -45& -2& 6& -10& -1& -1& 3& 1& -1& 2& -2& -2& 10& -30& -30& 1& -3& 0& 0& 0& 0& 30& 90\\ \hline
$\chi_{58}$&1890& -30& 0& 0& 2& -6& -6& 26& 3& -45& 2& -2& -2& -1& -1& 3& -1& 1& 2& -2& -2& 2& 66& -30& -1& 3& 0& 0& 0& 0& -30& 90\\ \hline
$\chi_{59}$&2268& -36& 0& 0& -12& 28& 12& -36& 0& 0& 0& 0& 0& 0& 0& 0& 0& 0& 0& 0& 0& 12& -36& -36& 0& 0& 0& 0& 0& 0& 36& -36\\ \hline
$\chi_{60}$&2520& -40& 0& 0& 8& -24& 8& -24& -1& 15& 4& -4& -4& 3& -1& -1& 1& -1& -4& 4& 4& -8& 24& 24& 1& -3& -6& 2& 0& 0& 40& 200\\ \hline
$\chi_{61}$&2520& -40& 0& 0& 8& 8& -24& 8& -1& 15& -4& 4& 4& -1& -1& 3& -1& 1& 4& -4& -4& -8& 24& 24& 1& -3& -6& 2& 0& 0& 40& -280\\ \hline
$\chi_{62}$&2520& -40& 0& 0& 8& -8& -8& -8& 2& -30& 0& 0& 0& -2& 2& -2& 0& 0& 0& 0& 0& -8& -72& 56& -2& 6& 3& -1& -1& 3& 40& -40\\ \hline
$\chi_{63}$&2520& -40& 0& 0& 8& -8& -8& -8& 2& -30& 0& 0& 0& -2& 2& -2& 0& 0& 0& 0& 0& -8& 120& -8& -2& 6& 3& -1& 1& -3& 40& -40\\ \hline
$\chi_{64}$&2835& -45& 0& 0& 3& -25& 15& -9& 0& 0& -3& 1& 9& 0& 0& 0& 0& 0& -3& 5& -3& 3& 27& -21& 0& 0& 0& 0& 0& 0& -45& -225\\ \hline
$\chi_{65}$&2835& -45& 0& 0& 3& 11& -21& 27& 0& 0& 3& -1& -9& 0& 0& 0& 0& 0& 3& -1& -9& 3& -45& 3& 0& 0& 0& 0& 0& 0& -45& 315\\ \hline
$\chi_{66}$&2835& -45& 0& 0& 3& 15& -9& -33& 0& 0& -3& 5& -3& 0& 0& 0& 0& 0& -3& 1& 9& 3& 27& -21& 0& 0& 0& 0& 0& 0& -45& 135\\ \hline
$\chi_{67}$&2835& -45& 0& 0& 3& 3& 3& -45& 0& 0& 3& -5& 3& 0& 0& 0& 0& 0& 3& -5& 3& 3& -45& 3& 0& 0& 0& 0& 0& 0& -45& -45\\ \hline
\end{tabular}\\[.2cm]

\setlength\tabcolsep{1pt}
\footnotesize
\begin{tabular}{*{40}{|c}}
\hline
&33&34&35&36&37&38&39&40&41&42&43&44&45&46&47&48&49&50&51&52&53&54&55&56&57&58&59&60&61&62&63&64&65&66&67 \\
\hline\hline
$\chi_{43}$&14& -1& 1& 6& -10& -2& -2& 1& -1& -6& 2& 0& 0& 0& 0& 0& 0& 0& 0& 2& 2& -6& 0& 0& -2& -2& 2& 2& 0& -1& -1& 7& 1& -1& 0\\ \hline
$\chi_{44}$&-18& 1& -1& -2& -10& 2& 10& 1& -1& -6& 2& 0& 0& 0& 0& 0& 0& 0& 0& 2& -2& 2& 0& 0& 2& -2& -2& 2& 0& 1& -3& 5& -1& 1& 0\\ \hline
$\chi_{45}$&-18& 1& -1& 6& -2& -6& 2& -1& 1& -6& 2& 0& 0& 0& 0& 0& 0& 0& 0& 2& -2& 2& 0& 0& -2& 2& 2& -2& 0& 1& -3& 5& -1& 1& 0\\ \hline
$\chi_{46}$&6& -1& 1& -2& -2& -2& 14& -1& 1& -6& 2& 0& 0& 0& 0& 0& 0& 0& 0& -6& 2& 2& 0& 0& -2& -2& 2& 2& 0& -1& 3& -5& -1& 1& 0\\ \hline
$\chi_{47}$&1& 0& 0& 1& 9& 1& -15& 0& 0& -3& 1& 0& 0& 0& 0& 0& 0& -1& 1& -3& 1& 1& 1& -1& -3& 1& 1& 1& 0& -3& 1& 9& -1& 1& 0\\ \hline
$\chi_{48}$&25& 0& 0& 1& 9& -7& 9& 0& 0& 9& -3& 0& 0& 0& 0& 0& 0& 1& -1& 1& -3& 5& -1& 1& 1& -3& 1& 1& 0& -3& 1& 9& -1& 1& 0\\ \hline
$\chi_{49}$&29& 0& 0& -7& 9& 5& -3& 0& 0& -3& 1& 0& 0& 0& 0& 0& 0& 1& -1& 1& 1& -3& 1& -1& -1& -1& 3& -1& 0& -3& 5& -3& 1& -1& 0\\ \hline
$\chi_{50}$&5& 0& 0& 1& -15& 5& -3& 0& 0& 9& -3& 0& 0& 0& 0& 0& 0& -1& 1& 5& -3& 1& -1& 1& -1& -1& -1& 3& 0& -3& 5& -3& 1& -1& 0\\ \hline
$\chi_{51}$&16& 2& -2& 0& 0& 0& 0& 0& 0& 0& 0& 0& 0& -1& 3& 1& -1& 0& 0& 0& 0& 0& 0& 0& 0& 0& 0& 0& 0& 2& -2& -2& 0& 0& 0\\ \hline
$\chi_{52}$&-20& 2& -2& -4& 12& 4& -12& 0& 0& 12& -4& 0& 0& 0& 0& 0& 0& 0& 0& -4& 4& -4& 0& 0& 0& 0& 0& 0& 0& 2& -2& -2& 0& 0& 0\\ \hline
$\chi_{53}$&-40& 0& 0& 0& 0& 0& 0& 0& 0& 0& 0& 0& 0& 1& -3& 1& -1& 0& 0& 0& 0& 0& 0& 0& 0& 0& 0& 0& 0& 3& -1& -9& -1& 1& 0\\ \hline
$\chi_{54}$&-8& 0& 0& 0& 0& 0& 0& 0& 0& 0& 0& 0& 0& 1& -3& 1& -1& 0& 0& 0& 0& 0& 0& 0& 0& 0& 0& 0& 0& 3& -5& 3& 1& -1& 0\\ \hline
$\chi_{55}$&-22& 0& 0& 2& -6& -2& 6& 0& 0& 6& -2& 0& 0& 0& 0& 0& 0& 0& 0& -2& 2& -2& 0& 0& 2& -2& 2& -2& 0& -3& 5& -3& -1& 1& 0\\ \hline
$\chi_{56}$&34& 0& 0& -6& -6& 10& -6& 0& 0& 6& -2& 0& 0& 0& 0& 0& 0& 0& 0& 6& -2& -2& 0& 0& -2& -2& 2& 2& 0& 3& -5& 3& -1& 1& 0\\ \hline
$\chi_{57}$&-38& 0& 0& 2& -6& -2& 6& 0& 0& 6& -2& 0& 0& 0& 0& 0& 0& 0& 0& -2& 2& -2& 0& 0& -2& 2& -2& 2& 0& -3& 1& 9& 1& -1& 0\\ \hline
$\chi_{58}$&26& 0& 0& 2& 18& -6& -6& 0& 0& 6& -2& 0& 0& 0& 0& 0& 0& 0& 0& -2& -2& 6& 0& 0& -2& -2& 2& 2& 0& 3& -1& -9& 1& -1& 0\\ \hline
$\chi_{59}$&-36& 0& 0& -4& 12& 4& -12& 0& 0& -12& 4& 0& 0& -1& 3& -1& 1& 0& 0& 4& -4& 4& 0& 0& 0& 0& 0& 0& 0& 0& 0& 0& 0& 0& 0\\ \hline
$\chi_{60}$&-56& -2& 2& 0& 0& 0& 0& 0& 0& 0& 0& 0& 0& 0& 0& 0& 0& 0& 0& 0& 0& 0& 0& 0& 0& 0& 0& 0& 0& 1& 1& -7& 1& -1& 0\\ \hline
$\chi_{61}$&-24& -2& 2& 0& 0& 0& 0& 0& 0& 0& 0& 0& 0& 0& 0& 0& 0& 0& 0& 0& 0& 0& 0& 0& 0& 0& 0& 0& 0& 1& -3& 5& -1& 1& 0\\ \hline
$\chi_{62}$&-40& 1& -1& 8& 8& -8& -8& 1& -1& 0& 0& 0& 0& 0& 0& 0& 0& 0& 0& 0& 0& 0& 0& 0& 0& 0& 0& 0& 0& -2& 2& 2& 0& 0& 0\\ \hline
$\chi_{63}$&-40& 1& -1& -8& -8& 8& 8& -1& 1& 0& 0& 0& 0& 0& 0& 0& 0& 0& 0& 0& 0& 0& 0& 0& 0& 0& 0& 0& 0& -2& 2& 2& 0& 0& 0\\ \hline
$\chi_{64}$&63& 0& 0& 3& 3& -1& -9& 0& 0& -9& 3& 0& 0& 0& 0& 0& 0& -1& 1& -5& 3& -1& -1& 1& 1& 1& -3& 1& 0& 0& 0& 0& 0& 0& 0\\ \hline
$\chi_{65}$&27& 0& 0& -5& 3& 3& 3& 0& 0& -9& 3& 0& 0& 0& 0& 0& 0& 1& -1& -1& 3& -5& -1& 1& 3& -1& -1& -1& 0& 0& 0& 0& 0& 0& 0\\ \hline
$\chi_{66}$&39& 0& 0& 11& -21& -1& -9& 0& 0& 3& -1& 0& 0& 0& 0& 0& 0& 1& -1& -1& -1& 3& 1& -1& 1& 1& 1& -3& 0& 0& 0& 0& 0& 0& 0\\ \hline
$\chi_{67}$&51& 0& 0& -5& 3& -5& 27& 0& 0& 3& -1& 0& 0& 0& 0& 0& 0& -1& 1& 3& -1& -1& 1& -1& -1& 3& -1& -1& 0& 0& 0& 0& 0& 0& 0\\ \hline\end{tabular}\\[.2cm]


\begin{thebibliography}{qcomp}

\bibitem[Ap]{Ap}
D. M. Appleby, SIC-POVMs and the extended Clifford group,
\href{https://doi.org/10.1063/1.1896384} 
{{\em J. Math. Phys.} {\bf 46}, 052107 (2005)}.


\bibitem[BOZ]{BOZ}
E. Bannai, M. Oura, and D. Zhao,
The complex conjugate invariants of Clifford groups,
\href{https://doi.org/10.1007/s10623-020-00819-7} 
{{\em Designs, Codes and Cryptography},
{\bf 89}(2) pp. 341–350, 2021}.


\bibitem[BRW]{BRW}
B. Bolt, T. Room, and G. Wall, 
On the Clifford collineation, transform and similarity groups. II.,
\href{https://doi.org/10.1017/S1446788700026380} 
{{\em J. Austral. Math. Soc.} {\bf 2}, pp. 80–96 (1961)}.


\bibitem[CTV]{CTV}
E. Campbell,  B. Terhal, and C. Vuillot, 
Roads towards fault-tolerant universal quantum computation, 
\href{https://doi.org/10.1038/nature23460} 
{{\em Nature} {\bf 549}, 172 (2017)}.




\bibitem[F+7]{F+7}
K. Fisher, A. Broadbent, L. Shalm, Z. Yan, J. Lavoie, R. Prevedel, T. Jennewein, and K. J. Resch, 
Quantum computing on encrypted data,
\href{https://doi.org/10.1038/ncomms4074} 
{{\em Nat. Commun.}, {\bf 5}, no. 1, p. 3074, 2014}.


\bibitem[FH]{FH}
W. Fulton and J. Harris,
{\em Representation theory: A first course}, 
\href{https://doi.org/10.1007/978-1-4612-0979-9}
{Springer, New York, 2013, Vol. 129}.




\bibitem[GAP]{GAP}
  The GAP~Group, \emph{GAP -- Groups, Algorithms, and Programming, 
  Version 4.12.2}; 
  2022,\\
  \url{https://www.gap-system.org}.


\bibitem[Go1]{Go1}
D. Gottesman, 
The Heisenberg Representation of Quantum Computers,
\href{https://doi.org/10.48550/arXiv.quant-ph/9807006}
{{\em Proc. of the 22nd Int. Colloq. on Group Theoretical Methods in Physics}}, eds. S. P. Corney, R. Delbourgo, and P. D. Jarvis
pp 32–43,  (Cambridge, MA, International Press, 1999). 

\bibitem[Go2]{Go2}
D. Gottesman, 
An introduction to quantum error correction and fault-tolerant quantum computation,
{\em Quantum Information Science and its Contributions to Mathematics}, 
\href{https://doi.org/10.48550/arXiv.0904.2557}
{Proc. Symp. in Applied Mathematics, 
{\bf 68} pp 13–58}.




\bibitem[GS]{GS}
D. Grier, and L. Schaeffer,
The Classification of Clifford Gates over Qubits,
\href{https://doi.org/10.22331/q-2022-06-13-734}
{{\em Quantum} {\bf 6}, 734 (2022)}.




\bibitem[Gr]{Gr}
D. Gross, 
Hudson’s theorem for finite-dimensional quantum systems. 
\href{https://doi.org/10.1063/1.2393152}
{{\em J. Math. Phys.} {\bf 47}(12), 122107 (2006)}. 


\bibitem[GAE]{GAE}
D. Gross, K. Audenaert, and J. Eisert, 
Evenly distributed unitaries: On the structure of unitary designs,
\href{https://doi.org/10.1063/1.2716992}
{{\em J. Math. Phys.} 48, 052104 (2007)}.


\bibitem[GNW]{GNW}
D. Gross, S. Nezami, and M. Walter, 
Schur-Weyl duality for the Clifford group with applications: Property testing, a robust Hudson theorem, and de Finetti representations,
\href{https://doi.org/10.1007/s00220-021-04118-7} 
{{\em Commun. Math. Phys.} {\bf 385}, 1325 (2021)}.






\bibitem[He]{He}
M. Heinrich, 
\href{https://quantumcomputing.stackexchange.com/questions/26150}{StackExchange discussion} on 
``Is the Clifford group perfect (equals its own commutator subgroup)?'' (2022-04-26)



\bibitem[HWFW]{HWFW}
J. Helsen, J. Wallman, S. Flammia, and S. Wehner,
Multiqubit randomized benchmarking using few samples,
\href{https://doi.org/10.1103/PhysRevA.100.032304}
{{\em Phys. Rev. A} {\bf 100}, 032304 (2019)}.






\bibitem[HWW]{HWW}
J. Helsen, J. Wallman, and S. Wehner, 
Representations of the multi-qubit Clifford group,
\href{https://doi.org/10.1063/1.4997688} 
{{\em J. Math. Phys.} {\bf 59}, 072201 (2018)}. 


\bibitem[JL]{JL}
G. James, and M. Liebeck,
{\em Representations and Characters of Groups},
\href{https://doi.org/10.1017/CBO9780511814532} 
{2nd edition, 2001, Cambridge University Press}.




\bibitem[KMP]{KMP}
C. Keeler, W. Munizzi, and J. Pollackb,
Clifford Orbits from Cayley Graph Quotients,
\href{https://doi.org/10.48550/arXiv.2306.01043}
{arXvi:2306.01043}




 
\bibitem[LC]{LC} 
C.-H. Lin and Y.-Y. Chen, 
HyperQUEEN: Hyperspectral Quantum Deep Network For Image Restoration, 
\href{https://doi.org/10.1109/TGRS.2023.3274355}
{{\em IEEE Transactions on Geoscience and Remote Sensing}, 
{\bf 61}, pp. 1-20, 2023}. 



\bibitem[Ma]{Ma}
K. Mastel,
The Clifford theory of the $n$-qubit Clifford group,
\href{https://doi.org/10.48550/arXiv.2307.05810}
{arXiv:2307.05810}



\bibitem[McBS+]{McBS+} 
J.R. McClean, S. Boixo, V.N. Smelyanskiy, et al., 
Barren plateaus in quantum neural network training landscapes,
\href{https://doi.org/10.1038/s41467-018-07090-4}
{{\em Nat. Commun.} {\bf 9}, 4812 (2018)}. 





\bibitem[MMG1]{MMG1}
F. Montealegre-Mora, and D. Gross, 
Rank-deficient representations in the theta correspondence over finite fields arise from quantum codes,
\href{https://doi.org/10.1090/ert/563}
{{\em Representation Theory} {\bf 25} pp.193--223, 2021}.


\bibitem[MMG2]{MMG2}
F. Montealegre-Mora, and D. Gross, 
Duality theory for Clifford tensor powers,
\href{https://doi.org/10.48550/arXiv.2208.01688}
{arXiv:2208.01688}.




\bibitem[NWD]{NWD}
P. Niemann, R. Wille, and R. Drechsler,
Efficient synthesis of quantum circuits implementing Clifford group operations,
\href{https://doi.org/10.1109/ASPDAC.2014.6742938} 
{{\em Asia and South Pacific Design Automation Conf.}, 
pp. 483-488, 2014}.



\bibitem[Oz]{Oz}
M. Ozols, 
Clifford group, 
Essays at University of Waterloo, Spring (2008).



\bibitem[Sa]{Sa}
J. Saied, 
{\em private communication} and discussion on
\href{https://scirate.com/arxiv/2309.14850}
{https://scirate.com/arxiv/2309.14850}


\bibitem[Sag]{Sag} 
B. Sagan,
{\em The Symmetric Group: Representations, Combinatorial Algorithms, and Symmetric Functions},
\href{https://doi.org/10.1007/978-1-4757-6804-6}
{2nd edition, GTM {\bf 203}, Springer-Verlag, New York, 2001}.


\bibitem[Sel]{Sel}
P. Selinger, 
Generators and relations for $n$-qubit Clifford operators, 
\href{https://doi.org/10.2168/LMCS-11(2:10)2015}
{{\em Logical Methods in Computer Science} {\bf 11} (2:10), pp.1--17 (2015)}.



\bibitem[SH]{SH}
T. Singal, and M.-H. Hsieh,
Approximate 3-designs and partial decomposition of the Clifford group representation using transvections,
\href{https://doi.org/10.48550/arXiv.2111.13678}
{arXiv.2111.13678}





\bibitem[WF]{WF}
J. Wallman, and S. Flammia, 
Randomized benchmarking with confidence,
\href{https://doi.org/10.1088/1367-2630/16/10/103032}
{{\em New Journal of Physics}, {\bf 16}(10):103032, 2014}. 


\bibitem[We]{We}
Z. Webb, 
The Clifford group forms a unitary 3-design,
\href{https://doi.org/10.26421/QIC16.15-16-8}
{{\em Quantum Information \& Computation}, 
{\bf 16}:1379, 2016}.


\bibitem[Zh]{Zh}
H. Zhu, 
Multiqubit Clifford groups are unitary 3-designs, 
\href{https://doi.org/10.1103/PhysRevA.96.062336}
{{\em Phys. Rev. A} {\bf 96}, 062336 (2017)}.

\bibitem[ZKGG]{ZKGG}
H. Zhu, R. Kueng, M. Grassl, and D. Gross, 
The Clifford group fails gracefully to be a unitary 4-design, 
\href{https://doi.org/10.48550/arXiv.1609.08172}
{arXiv:1609.08172} 

\end{thebibliography}
\end{document}